\documentclass[11pt]{amsart}
\usepackage{amssymb,amsmath,amscd}
\usepackage{amsthm}
\usepackage{amsbsy}
\usepackage{enumerate}
\usepackage{pstricks}
\usepackage{pst-plot}
\usepackage[all]{xy}

\def\ord{\text{ord}}
    \def\nz{\hbox{\text{0}{ \raise 3pt \hbox{\kern -10pt
    \vrule width 8pt height 0.4 pt }\kern 2pt}}}

\begin{document}
\title{EXISTENCE OF DICRITICAL DIVISORS REVISITED}
\author{By Shreeram S. Abhyankar and William J. Heinzer}
\address{Mathematics Department, Purdue University, West Lafayette, IN 47907, USA.}
\email{ram@cs.purdue.edu, heinzer@math.purdue.edu}
\begin{abstract}
We characterize the dicriticals of special pencils.
We also initiate higher dimensional dicritical theory.

2000 Mathematics Subject Classification: Primary 14A05.

Key Words: dicritical, special pencils, normal varieties.
\end{abstract}
\date{}
\maketitle

\markboth{}{}

\baselineskip 17pt

{{\bf Section 1: Introduction.}} 

The analytical (topological) theory of dicritical divisors was developed in
\cite{Art,EiN,Fou,LeW,MaM}. It was algebracized in \cite{Ab8,Ab9}. The algebraic
theory was furthered in \cite{Ab10,Ab11,AbH,AbL}.  In this paper we shall make 
further progress in this theory. In particular, in Theorem (8.2) we shall prove 
a converse of Proposition (3.5) of \cite{Ab11} characterizing the dicritical set 
of a special pencil on a nonsingular surface; see Section 8 for the statement of 
this and other related results. In Section 9 we shall initiate the dicritical 
theory of higher dimensional normal varieties, which may be algebraic or 
arithmetical, and we shall indicate how this could be used in attacking the 
higher dimensional jacobian conjecture.

We shall use the notation and terminology introduced in \cite{Ab3} to \cite{Ab9}
and more specifically in \cite{Ab10,Ab11}. Relevant
background material can be found in \cite{Ab1,Ab2,Nag,NoR,Zar}.

It may be pointed out that the theory of graded rings, say as found on pages
206-215, 236-241, 272-277, 399-408, and  585-587 of \cite{Ab4}, is a cornerstone 
of this paper.  The basis of that theory is the idea of collecting terms of 
like degree in a polynomial coming from classical algebra. This idea is used
in geometry in the aphorism which says that the factors of the highest degree
terms of a bivariate polynomial give points at infinity while the factors of
its lowest degree term give tangents at the origin.
Another cornerstone of this paper consists of the theories of blowing up and
Veronese embedding as developed on pages 7-45, 155-192, 262-283
of \cite{Ab3} and reproduced on pages 146-161, 534-552, 553-577 of \cite{Ab4}.

In Sections 2 to 5 we shall review some notation to be used frequently.
In Section 6, which is the heart of the paper, we study natural extensions
of valuations, which are sometimes called Gauss extensions. Behind all this
are Rees rings and their suitable homomorphic images which we call form rings
and which are sometimes called fiber rings.
In Section 7 we make connections with extended Rees rings.

\centerline{}

{{\bf Section 2: Quasilocal Rings.}} 

Recall that: 
A ring is commutative with $1$.
A quasilocal ring is a ring $S$ with a unique maximal ideal $M(S)$; by
$$
H_S:S\to H(S)=S/M(S)
$$
we denote the residue class epimorphism. 
A quasilocal ring $S$ dominates a 
quasilocal ring $T$ means $T$ is a subring of $S$ with $M(T)=T\cap M(S)$,
and then we say that $S$ is residually rational (resp: residually algebraic,
residually transcendental, residually simple transcendental, residually almost
simple transcendental, ...) over $T$ if $H(S)=H_S(T)$
(resp: $H(S)$ is algebraic over $H_S(T)$,
$H(S)$ is transcendental over $H_S(T)$,
$H(S)$ is simple transcendental over $H_S(T)$,
$H(S)$ is almost simple transcendental over $H_S(T)$, ...). Note that a field
$k^*$ is simple transcendental over a subfield $k$ means $k^*=k(t)$ for some
element $t$ which is transcendental over $k$, and a field $k^*$ is 
{\bf almost simple transcendental} over a subfield $k$ means $k^*$ is simple
transcendental over a finite algebraic field extension of $k$.
A local ring is a noetherian quasilocal ring.

As usual $\mathbb N$ (resp: $\mathbb N_+$) denotes the set of all nonnegative 
(resp positive) integers. The set of all nonzero elements in a ring $A$ is denoted
by $A^\times$. 

 For any set $U$ of quasilocal domains and any $i\in\mathbb N$,
$U_i$ denotes the set of all $i$-dimensional members of  $U$.

See pages 115 and 127 of \cite{Ab4} for the definitions of:
dim$(R)$, spec$(R)$, mspec$(R)$, ht$_RJ$, dpt$_RJ$, vspec$_RJ$, mvspec$_RJ$, and
nvspec$_RJ$, for a ring $R$ and and an ideal $J$ in it. Note that
nvspec$_RJ$ stand for the minimal spectrum of $J$ in $R$, and its members are 
called the minimal primes of $J$ in $R$. For the definition of
vector space dimension 
$$
[L:K]=\dim_KL
$$ 
see page 9 of \cite{Ab4}.

See the third paragraph of Section 3 of \cite{Ab10} for the definition of
coefficient set, coefficient ring, and coefficient field of a quasilocal
ring.

\centerline{}

{{\bf Section 3: Modelic Blowups.}} 

Referring to pages 145-161 of \cite{Ab4}, for 
the foundations of models, recall that: 

$\mathfrak V(A)=\{A_P:P\in\text{spec}(A)\}=$ the modelic spec of a domain $A$.

$S^{\mathfrak N}=$ the set of all members of $\mathfrak V(\overline S)$ which
dominate $S$ where $\overline S$ is the integral closure of a quasilocal domain
$S$ in its quotient field QF$(S)$.

$\mathfrak W(A,J)=\cup_{0\ne x\in J}\mathfrak V(A[Jx^{-1}])=$ the modelic 
blowup of $A$ at a nonzero ideal $J$ in a domain $A$; note that
$Jx^{-1}=\{y/x:y\in J\}$; see pages 152-160 of \cite{Ab4}.

$\mathfrak W(S,J)^\Delta_i=$ the set of all $i$-dimensional members of
$\mathfrak W(S,J)$ which dominates $S$ where $J$ is a nozero ideal in a quasilocal
domain $S$ and $i\in\mathbb N$.

$\mathfrak D(S,J)=(\mathfrak W(S,J)_1^\Delta)^{\mathfrak N}=$ the 
{\bf dicritical set} of a nonzero ideal $J$ in a quasilocal domain $S$; members of
this set are called {\bf dicritical divisors} of $J$ in $S$. 

$\overline C(A)=$ the set of all nonzero complete ideals in a normal domain $A$.
See the second paragraph after (2.3) of \cite{Ab11} for the definitions of:
normal domain, simple ideal, valuation ideal, complete ideal, normal ideal, and
completion of an ideal.

$C(S)=$ the set of all $M(S)$-primary simple complete ideals in a normal
local domain $S$.

ord$_Sa
=\text{max}\{d\in\mathbb N:a\in M(S)^d\text{ or }a\subset M(S)^d\}$ 
where $S$ is a local ring and $a\in S$ or $a\subset S$, with the
understanding that if $a=0$ or $a\subset\{0\}$ then
ord$_Sa=\infty$.

ord$_S(\alpha/\beta)=(\ord_S\alpha)-(\ord_S\beta)$ for $\alpha\ne 0\ne\beta$
in a regular local domain $S$.

DVR = Discrete Valuation Ring = one dimensional regular local domain.

$o(S)=$ the natural DVR of a positive dimensional regular local domain $S$, 
i.e., the DVR $V$ with QF$(V)=L=\text{QF}(S)$ such that 
ord$_V(a)=\ord_S(a)$ for all $a\in L$. 

$V(a)=v(a)$ for all 
$a\in L$ where $v:L\to G\cup\{\infty\}$ is a valuation of a field $L$ and
$V$ is its valuation ring. Note that $V(a)=\infty$ or
$V(a)\in v(L^\times)$ according as $a=0$ or $a\ne 0$; page see 41 of \cite{Ab4}.
Now let $I\subset S$ where $S$ is a noetherian subring of $L$ with $S\subset V$.
If $I\subset\{0\}$ then we let $V(I)=\infty$.
If $I\not\subset\{0\}$ then we let $V(I)$ be the unique element of $v(L^\times)$
such that for some $0\ne x\in I$ we have $V(I)=v(x)\le v(y)$ for all $y\in I$.
Note that the noetherianness guarantees the existence of $x$.
Also note that the equation $V(I)=V(x)$ is equivalent to the equation $IV=xV$.

\centerline{}

{{\bf Section 4: Graded Rings.}} Referring to pages 206-215 of \cite{Ab4}, for the
foundations of graded rings, recall that: 

$E_R(I)=R[IZ]=$ the Rees ring of an ideal $I$ in a nonnull ring $R$ relative 
to $R$ with variable $Z$. Note that $R[Z]$ is the univariate polynomial ring as a 
naturally graded homogeneous ring with $R[Z]_n$ = the set of all homogeneous
polynomials of degree $n$ including the zero polynomial, and $n$ varying over
$\mathbb N$. Now $E_R(I)$ is a graded subring of $R[Z]$ with
$E_R(I)_n=\{gZ^n:g\in I^n\}$, 
and every $f\in E_R(I)$ can uniquely be written as a finite sum
$$
f=\sum_{n\in\mathbb N}f_nZ^n\;\;\text{ with }\;\; f_n\in I^n.
\leqno(4.1)
$$
Details are in the third paragraph after (2.3) of \cite{Ab11} where you can also
find the the definitions of:
an element or subset of a nonnull ring $S$ to be integral over an ideal $J$ in a
subring of $S$, the integral closure of $J$ in $S$, and the reduction of an ideal.

$F_R(I)=E_R(I)/ME_R(I)=$ the form ring of an ideal $I\subset M=M(R)$ in a local 
ring $R$ relative to $R$ with variable $Z$. Note that $ME_R(I)$ is a homogeneous 
ideal in $E_R(I)$ and hence $F_R(I)$ is a naturally graded homogeneous ring over
the field $R/M$, and for its homogeneous 
$n$-th component $F_R(I)_n$ we have a canonical $R$-epimorphism
$$ 
\mu_n:I^n\to F_R(I)_n
\leqno(4.2)
$$
with kernel $MI^n$. Details are in the third paragraph after (2.3) of \cite{Ab11},
where we slightly generalized the matter by letting $I$ to be an ideal in a
nonnull ring $R$ with $I\subset M=$ a nonunit ideal in $R$, and denoting the
form ring by $F_{(R,M)}(I)$.
Note that $F_{(R,M)}(I)$ is isomorphic as a graded ring to the
associated graded ring grad$(R,I,M)$ of Definition (D3) on page 586 of \cite{Ab4}.

\centerline{}

{{\bf Section 5: Quadratic Transformations.}} 
Let $R$ be a two dimensional regular local domain with quotient field QF$(R)=L$. 
Recall that: 

$D(R)^\Delta=$ the set of all DVRs $V$ with quotient field $L$ such that $V$ 
dominates $R$ and is residually transcendental over $R$; 
by \cite{Ab1,Ab2} it follows that then $H(V)$ is almost simple transcendental over 
$H_V(R)$. Members of $D(R)^\Delta$ are called prime divisors of $R$.

QDT = Quadratic Transformation or Quadratic Transform.

$Q_j(R)=$ the set of all two dimensional QDTs of $R$.

$Q(R)=\coprod_{j\in\mathbb N}Q_j(R)=$ the set of all two dimensional regular
local domains which birationally dominate $R$, i.e., whose quotient field is
$L$ and which dominate $R$; proof of the second equality in \cite{Ab1}.

$o_R:Q(R)\to D(R)^\Delta$ is the bijective map given by $S\mapsto o(S)$.
The QDT sequence of $R$ along $o(S)$ is the finite sequence $(R_j)_{0\le j\le\nu}$
with $R_0=R$ and $R_\nu=S$ such that $R_{j+1}\in Q_1(R_j)$ for $0\le j<\nu$. 
A finite QDT sequence of $R$ means the QDT sequence of $R$ along some 
$V\in D(R)^\Delta$.  Proofs in \cite{Ab1}.

$\zeta_R:D(R)^\Delta\to C(R)$ is the bijective map given in Appendix 5 of
\cite{Zar} which we call the Zariski map. Details are in Section 2 of \cite{Ab11}.

$a_R(z)$ (resp: $b_R(z),J_R(z),I_R(z))=$ the numerator (resp: denominator, 
first asociated, second associated) ideal of $z$ in $R$.
Details are in the fourth paragraph after (2.3) of \cite{Ab11}.

$\mathfrak D(R,z)^\sharp\subset\mathfrak D(R,z)^\flat\subset
\mathfrak D(R,z)=(\mathfrak W(R,J_R(z))^\Delta_1)^{\mathfrak N}$
where $\mathfrak D(R,z)^\sharp$ (resp: $\mathfrak D(R,z)^\flat$,
$\mathfrak D(R,z)$) is the set of all {\bf sharp dicritical divisors}
(resp: {\bf flat dicritical divisors, dicritical divisors}) of $z$ 
in $R$, i.e., 
the set of those DVRs
$V\in\ D(R)^\Delta$ at which the element $z\in L^\times$ 
is a residual transcendental generator (resp: residually a polynomial, 
residually transcendental) over $R$.
We call $\mathfrak D(R,z)^\sharp$ (resp: $\mathfrak D(R,z)^\flat$,
$\mathfrak D(R,z)$) the {\bf sharp dicritical set} (resp: 
the {\bf flat dicritical set}, the {\bf dicritical set}) of $z$ in $R$.
See the material starting with the second display in Section 2 of \cite{Ab10}.

Geometrically speaking, we may visualize $R$ to be the local ring of a simple
point of an algebraic or arithmetical surface, and think of $z$ as a 
{\bf rational function} at that simple point which corresponds to the 
{\bf local pencil} of curves $a=ub$ at that point.  We say that $z$ generates a 
{\bf special pencil at} $R$ to mean that $b$ can be chosen so that $b=x^m$ for 
some $x\in M(R)\setminus M(R)^2$ and $m\in\mathbb N$, i.e., $zx^m\in R$ for some 
$x\in M(R)\setminus M(R)^2$ and $m\in\mathbb N$.  We say that $z$ generates a 
{\bf semispecial pencil at} $R$ to mean that $b$ can be chosen so that $b=x^my^n$ 
for some $x,y$ in $M(R)$ and $m,n$ in $\mathbb N$ with $M(R)=(x,y)R$, i.e., 
$zx^my^n\in R$ for some $x,y$ in $M(R)$ and $m,n$ in 
$\mathbb N$ with $M(R)=(x,y)R$. 
We say that $z$ generates a {\bf polynomial} or {\bf nonpolynomial} pencil in $R$
according as $\mathfrak D(R,z)=\mathfrak D(R,z)^\flat$ or
$\mathfrak D(R,z)\ne\mathfrak D(R,z)^\flat$.
We say that $z$ generates a {\bf generating} or {\bf nongenerating} pencil in $R$
according as $\mathfrak D(R,z)=\mathfrak D(R,z)^\sharp$ or
$\mathfrak D(R,z)\ne\mathfrak D(R,z)^\sharp$.
See the material starting with the second display in Section 2 of \cite{Ab10}.

 For a moment let $J$ be a nonzero ideal in $R$. We call $J$ a {\bf pencil} 
(in $R$) if $J=yJ_R(z)$ for some $y\in R^\times$ and $z\in L^\times$, and we note
that then $\mathfrak D(R,J)=\mathfrak D(R,z)$. If $J$ is a pencil with 
$J=yJ_R(z)$ then we let $\mathfrak D(R,J)^\sharp=\mathfrak D(R,z)^\sharp$ and 
$\mathfrak D(R,J)^\flat=\mathfrak D(R,z)^\flat$, and if $J$ is not a pencil then
we let $\mathfrak D(R,J)^\sharp=\mathfrak D(R,J)^\flat=\emptyset$.
Note that now we have 
$$
\mathfrak D(R,J)^\sharp\subset\mathfrak D(R,J)^\flat\subset
\mathfrak D(R,J)=(\mathfrak W(R,J)^\Delta_1)^{\mathfrak N}.
$$
We say that $J$ is a {\bf polynomial} or {\bf nonpolynomial} ideal in $R$
according as $J=yJ_R(z)$ for some $y\in R^\times$ and $z\in L^\times$ such that
$z$ generates a polynomial or nonpolynomial pencil in $R$. 
We say that $J$ is a {\bf generating} or {\bf nongenerating} ideal in $R$
according as $J=yJ_R(z)$ for some $y\in R^\times$ and $z\in L^\times$ such that
$z$ generates a generating or nongenerating pencil in $R$. 
In the last two sentences we may say {\bf pencil} instead of ideal.
Regardless of whether $J$ is a pencil or not, 
we say that $J$ is {\bf primary} to mean that the ideal $J$ is $M(R)$-primary.
We say that $J$ is {\bf special} (resp: {\bf semispecial}) {\bf at} $R$ if 
$J=yJ_R(z)$ for some $y\in R^\times$ and $z\in L^\times$ such that $z$ generates 
a special (resp: semispecial) pencil at $R$. 
We put
$$
\begin{cases}
\mathfrak B(R,J)^\sharp=\{o_R^{-1}(V):V\in\mathfrak D(R,J)^\sharp\}\\
\mathfrak B(R,J)^\flat=\{o_R^{-1}(V):V\in\mathfrak D(R,J)^\flat\}\\
\mathfrak B(R,J)=\{o_R^{-1}(V):V\in\mathfrak D(R,J)\}\\
\mathfrak Q(R,J)=\{T\in Q(R):(R,T)(J)\text{ is not principal}\}\\
\end{cases}
$$
and we note that $\mathfrak Q(R,J)$ is a finite set (for a proof see 
Proposition 2 on page 367 of \cite{Zar}) with
$$
\mathfrak B(R,J)^\sharp\subset\mathfrak B(R,J)^\flat
\subset\mathfrak B(R,J)\subset\mathfrak Q(R,J)\subset Q(R).
$$
We say that $J$ {\bf goes through} the members of $\mathfrak Q(R,J)$ but not
through the members of $Q(R)\setminus\mathfrak Q(R,J)$. 
See the middle of Section 2 of \cite{Ab10}.

\centerline{}

{{\bf Section 6: Natural Extensions of Valuations.}} Let $Y$ be an indeterminate 
over a field $L$.  Let $v:L\to G\cup\{\infty\}$ be a valuation of $L$. 
By (J4) to (J9) on pages 79-80 of \cite{Ab4} we get a unique valuation 
$w:L(Y)\to G\cup\{\infty\}$ of $L(Y)$ such that for all 
$\sum_{i\in\mathbb N} a_iY^i\in L[Y]$ 
with $a_i\in L$ we have
$$
w(\sum_{i\in\mathbb N} a_iY^i)=\text{min}\{v(a_i):i\in\mathbb N\}
\leqno(6.1)
$$ 
with the understanding that $w(0)=\infty$.
We call $w$ the {\bf natural extension} of $v$ to $L(Y)$; see \cite{Ab12}. We
rename the valuation rings $R_v$ and $R_w$ by putting
$$
V=R_v\quad\text{ and }\quad W=R_w
$$
and we call $W$ the {\bf natural extension} of $V$ to $L(Y)$.
Note that if $v$ is real discrete, i.e., if $G=v(L^\times)=\mathbb Z$, then 
$V$ and $W$ are DVRs.

\centerline{}

Let $R$ be a noetherian domain with quotient field $L$.
Let $I$ be a nonzero ideal in $R$.  Let 
$$
\text{$E=E_R(I)=R[IZ]=\sum_{n\in\mathbb N}E_n=$ the Rees ring of $I$.}
$$
Assume that $R\subset V$. We claim that then there exists a
unique valuation 
$$
w':L(Z)\to G\cup\{\infty\}
$$ 
such that, for all $0\ne f\in E$, in the notation of (4.1) we have
$$
w'(f)=\text{min}\{v(f_n)-nV(I):n\in\mathbb N\text{ with }f_n\ne 0\}.
\leqno(6.2)
$$ 
Namely, we can take $0\ne x\in I$ with $V(x)=V(I)$, and for any such $x$, upon
letting $Y=xZ$, by (6.1) and (6.2) we get 
$$
w'=w.
\leqno(6.3)
$$ 
So we call $w$ the $(R,I)$-{\bf extension} of $v$
to $L(Z)$ and we call $W$ the $(R,I)$-{\bf extension} of $V$ to $L(Z)$.

\centerline{}

Let $A$ be a noetherian domain with quotient field $L$, and let us consider the
meromorphic polynomial ring $B=A[Y,Y^{-1}]$. Also consider the multiplicative set
$M^*=\{1,Y,Y^2,\dots\}$ in the usual polynomial ring $B^*=A[Y]$.
By T(30.1) on page 233 of \cite{Ab4} we see that 
$$
\begin{cases}
\text{for any $P\in\text{spec}(A)$ we have:}\\
\quad\text{$PB^*\in\text{spec}(B^*)$ with $M^*\cap(PB^*)=\emptyset$}\\
\quad\text{ and }P=A\cap(PB^*)\text{ with }\text{ht}_AP=\text{ht}_{B^*}(PB^*).
\end{cases}
\leqno(6.4)
$$
Clearly $B$ equals the localization $B^*_{M^*}$ and hence by taking 
$(B^*,M^*)=(R,S)$ in (T12) on page 139 of \cite{Ab4} we get
$$
\begin{cases}
\text{for any $Q^*\in\text{spec}(B^*)\text{ with }M^*\cap Q^*=\emptyset$
we have:}\\
\quad Q^*B\in\text{spec}(B)\text{ with }M^*\cap(Q^*B)=\emptyset\\
\quad\text{ and }Q^*=B^*\cap(Q^*B)\text{ with }
\text{ht}_{B^*}Q^*=\text{ht}_{B}(Q^*B).
\end{cases}
\leqno(6.5)
$$
Taking $Q^*=PB^*$ in (6.5), by (6.4) and (6.5) we see that
$$
\begin{cases}
\text{for any $P\in\text{spec}(A)$ we have:}\\
\quad PB\in\text{spec}(B)\text{ with }M^*\cap(PB)=\emptyset\\
\quad\text{and }P=A\cap(PB)\text{ with } \text{ht}_AP=\text{ht}_{B}(PB).
\end{cases}
\leqno(6.6)
$$
Since every element of $L[Y,Y^{-1}]$ lands in $L[Y]$ after we multiply it by a 
high enough power of $Y$, by (6.1) it follows that for any
$\sum_{i\in\mathbb Z}a_iY^i\in L[Y,1/Y]^\times$ with $a_i\in L$ we have
$$
w(\sum_{i\in\mathbb Z} a_iY^i)=\text{min}\{v(a_i):i\in\mathbb Z\}.
\leqno(6.7)
$$ 

\centerline{}

Note that if $A\subset V$ then, upon letting $P=A\cap M(V)$ and $Q=B\cap M(W)$, 
we clearly have $P=A\cap Q\in\text{spec}(A)$ with $Q\in\text{spec}(B)$, and 
by (6.7) we get $B\subset W$ with $PB=Q$. Therefore by (6.6) it follows that
$$
\begin{cases}
\text{if $A\subset V$ then,}\\
\text{upon letting }P=A\cap M(V)\text{ and }Q=B\cap M(W),\text{ we have:}\\
\quad B\subset W\text{ with }P=A\cap Q\in\text{spec}(A)
\text{ and }M^*\cap Q=\emptyset\\
\quad\text{ and }Q=PB\in\text{spec}(B)\text{ with }
\text{ht}_AP=\text{ht}_{B}Q.
\end{cases}
\leqno(6.8)
$$

\centerline{}

Having separately dealt with the strands of the two subdomains $R$ and $A$
of $L$, let us now weave them together. 

\centerline{}

LEMMA (6.9).  Assume that $Y=xZ$ and $A=R[Ix^{-1}]$ where $0\ne x\in I$ is
such that $V(x)=V(I)$. Then 
$$
A\subset V\;\text{ and }\; E\subset B\subset W
\leqno(1)
$$
and
$$
\text{$M^*$ is a multiplicative set in $E$ with $E_{M^*}=B$}
\leqno(2)
$$ 
and upon letting
$$
P=A\cap M(V)\;\text{ and }\; Q=B\cap M(W)
\;\text{ and }\; P^*=E\cap M(W)
$$ 
we have
$$
P^*\in\text{spec}(E)
\leqno(3)
$$
and
$$
P=A\cap Q\in\text{spec}(A)\;\text{ with }\;Q\in\text{spec}(B)
\leqno(4)
$$
and
$$
M^*\cap Q=\emptyset
\leqno(5)
$$
and
$$
PB=Q=P^*B\;\text{ with }\;\text{ht}_AP=\text{ht}_BQ=\text{ht}_EP^*.
\leqno(6)
$$
 Finally, if $V\in(\mathfrak W(R,I)_1)^{\mathfrak N}$ then we have
$$
\text{ht$_AP=\text{ht}_BQ=\text{ht}_EP^*=1$.}
\leqno(7)
$$

\centerline{}

PROOF. Clearly $A\subset V$ and hence by (6.8) we get $B\subset W$.

Now $Y=xZ\in E_1\subset E$ and hence $M^*$ is a multiplicative set in $E$.
Every element of $E_1$ can be written as $yZ$ with $y\in I$, and we have
$yZ=(yx^{-1})Y$ with $yx^{-1}\in A$ and hence $E_1\subset A[Y]$.
Also $E_0=R\subset A$ and therefore $E\subset A[Y]$. Consequently
$E_{M^*}=E[Y^{-1}]\subset A[Y,Y^{-1}]=B$. 

To prove the reverse inclusion, note that $R\subset E$ 
and for every $y\in I$ we have
$(yx^{-1})Y=yZ\in E_1\subset E$. Consequently, for every $a\in A$ we have
$aY^m\in E$ for some $m\in\mathbb N$. Therefore 
$B=A[Y,Y^{-1}]\subset E[Y^{-1}]=E_{M^*}$. 

This proves (1) and (2). Now (3) and (4) are obvious.

By (6.8) we get (5). By (6.8) we also get
the first equalities in the two assertions
$$
\text{$PB=Q=Q^*B$\;\; and\;\; ht$_AP=\text{ht}_BQ=\text{ht}_EQ^*$}
$$ 
of (6), whereas the second 
equalities in these assertions follow from (2), (4), and (5), 
by invoking (T12) on page 139 of \cite{Ab4}. This proves (6).

If $V\in(\mathfrak W(R,I)_1)^{\mathfrak N}$ then clearly ht$_AP=1$, and hence
by (3) and (6) we get (7).

\centerline{}

{\it Remark $(6.9^\flat)$.} The above proof is not difficult but, 
because of the mixing of two 
strands, it is certainly subtle. Thus first we go up the ladder 
$A\subset B^*\subset B$ with prime ideals $P\subset Q^*\subset Q$
and then down the ladder $R\subset E\subset B$ with prime ideals
$R\cap M(V)\subset P^*\subset Q$.
This enables us to compare the ideal theories of the two seemingly uncomparable 
rings $A$ and $E$, neither of which is contained in the other.
This subtlety is accentuated in the proof of the following Lemma (6.11).
The subtlety of these proofs reminds me of the engraving which I had seen in Fine 
Hall of Princeton University Mathematics Department citing Einstein's quotation
``Raffiniert ist der Herr Gott, aber boshaft ist er nicht." 

\centerline{}

DEFINITION-OBSERVATION (6.10). Inspired by (6.9)(7), for any nonzero ideals 
$J\subset M$ in a domain $S$ and any $i\in\mathbb N$ we put

$\mathfrak W(S,J,M)^\Delta_i=$ the set of all $i$-dimensional members $T$ of
$\mathfrak W(S,J)$ such that $M\subset M(T)$.

We also put

$\mathfrak D(S,J,M)=(\mathfrak W(S,J,M)^{\Delta}_1)^{\mathfrak N}$ which we call
the {\bf dicritical set} of $(J,M)$ in $S$, and we call its members the 
{\bf dicritical divisors} of $(J,M)$ in $S$. 

We make the following observations concerning these concepts.

(I) If $N$ is an ideal in $S$ with $J\subset N$ such that rad$_SN=\text{rad}_SM$
then we have $\mathfrak W(S,J,N)^\Delta_i=\mathfrak W(S,J,M)^\Delta_i$ and
$\mathfrak D(S,J,N)=\mathfrak D(S,J,M)$.
If $M=S$ then $\mathfrak W(S,J,M)^\Delta_i=\emptyset$ and
$\mathfrak D(S,J,M)=\emptyset$.

(II) If $S$ is quasilocal and $M=M(S)$ then 
$\mathfrak W(S,J,M)^\Delta_i=\mathfrak W(S,J)^\Delta_i$ and
$\mathfrak D(S,J,M)=\mathfrak D(S,J)$.

(III) If $R$ is a noetherian domain and $I\subset M$ are nonzero nonunit ideals 
in $R$ then $\mathfrak D(R,I,M)$ is a 
finite set of DVRs. As in $(5.6)(\dagger^*)$ of \cite{Ab8},
this follows from (33.10) on page 118 of \cite{Nag} or
(33.2) on page 115 of \cite{Nag}. Upon letting $V_1,\dots,V_h$ be all the
distinct members of $\mathfrak D(R,I,M)$ and upon letting $W_1,\dots,W_h$ be their
respective $(R,I)$-extension to $L(Z)$, it follows that $W_1,\dots,W_h$ are
distinct DVRs with $E\subset W_j$ for $1\le j\le h$.
Note that if $M\subset\text{rad}_RI$ then clearly $h>0$.

\centerline{}

LEMMA (6.11).  Assume that $R$ is a normal noetherian domain and let $I$ be any 
nonzero nonunit normal ideal in $R$. Then $E$ is a normal noetherian domain
with $\mathfrak W(R,I)^{\mathfrak N}=\mathfrak W(R,I)$, 
$\mathfrak D(R,I,I)$ is a nonempty finite set of DVRs, and upon letting
$V_1,\dots,V_h$ be all the distinct members of $\mathfrak D(R,I,I)$, and upon 
letting $W_1,\dots,W_h$ be their respective $(R,I)$-extensions to $L(Z)$, we have 
that $W_1,\dots,W_h$ are distinct DVRs with $E\subset W_j$ for $1\le j\le h$. 

Moreover, upon letting 
$$
P^*_j=E\cap M(W_j)\;\;\text{ and }\;\;\overline P^*_j=E\cap(IW_j)\;\;
\text{ and }\;\;J_j=R\cap(IV_j)
$$ 
we have that $P^*_1,\dots,P^*_h$ are all the
distinct members of nvspec$_E(IE)$ and 
$$
\text{ht$_EP^*_j=1$\;\; with \;\;$E_{P^*_j}=W_j$\;\;
and\;\;$\overline P^*_j$ is $P^*_j$-primary.}
$$
for $i\le j\le h$.  Furthermore
$$
IE=\overline P^*_1\cap\dots\cap\overline P^*_h
$$ 
is the unique irredundant primary decomposition of $IE$ in $E$, and we have
$$
I=J_1\cap\dots\cap J_h\;\;\text{ with }\;\;
J_j=R\cap\overline P^*_j\;\;\text{ for }\;\;1\le j\le h.
$$

Now assume that $Y=xZ$ and $A=R[Ix^{-1}]$ where $0\ne x\in R$ is such that
$V_j(x)=V_j(I)$ for all $j$ in a nonempty subset $\Lambda$ of $\{1,\dots,h\}$.
Then $A$ and $B$ are normal noetherian domains with $A\subset V_j$ and 
$E\subset B\subset W_j$ for all $j\in\Lambda$. Moreover,
upon letting $P_j=A\cap M(V_j)$ with $\overline P_j=A\cap(IV_j)$ and
$Q_j=B\cap M(W_j)$ with $\overline Q_j=B\cap(IW_j)$, for all $j\ne j'$ in
$\Lambda$ we have 
$$
P_j\ne P_{j'}\;\;\text{ with }\;\;\overline P_j\ne\overline P_{j'}\;\;
\text{ and }\;\;Q_j\ne Q_{j'}\;\;\text{ with }\;\;\overline Q_j\ne\overline Q_{j'}.
$$
Furthermore, for all $j\in\Lambda$ we have that
$$
xA=IA\subset P_j\in\text{nvspec}_A(IA)\;\;\text{ and }\;\; 
xB=IB\subset Q_j\in\text{nvspec}_B(IB)
$$
and
$$
\text{ht$_AP_j=1$\;\; with \;\;$A_{P_j}=V_j$\;\;
and\;\;$\overline P_j$ is $P_j$-primary}
$$
and
$$
\text{ht$_BQ_j=1$\;\; with \;\;$B_{Q_j}=W_j$\;\;
and\;\;$\overline Q_j$ is $Q_j$-primary}
$$
and
$$
P_jB=Q_j=P^*_jB\;\;\text{ with }\;\;P_j=A\cap Q_j
\;\;\text{ and }\;\;P^*_j=E\cap Q_j
$$
and
$$
\overline P_jB=\overline Q_j=\overline P^*_jB
\;\;\text{ with }\;\;\overline P_j=A\cap \overline Q_j
\;\;\text{ and }\;\;\overline P^*_j=E\cap \overline Q_j.
$$
 Finally, if for all $l\in\{1,\dots,h\}\setminus\Lambda$ we have 
$V_l(x)\ne V_l(I)$, then
$$
IA=\cap_{j\in\Lambda}\overline P_j\;\;\text{ and }\;\;
IB=\cap_{j\in\Lambda}\overline Q_j
$$ 
are the unique irredundant primary decompositions of $IA$ and $IB$ in $A$ and $B$ 
respectively, and we have
$$
\text{nvspec}_A(IA)=\{P_j:j\in\Lambda\}
\;\;\text{ and }\;\;
\text{nvspec}_B(IB)=\{Q_j:j\in\Lambda\}.
$$ 

\centerline{}

PROOF. By (8.1)(VI) of \cite{AbH} and the above Lemma (6.9)(III) 
we see that: $E$ is a normal noetherian domain
with $\mathfrak W(R,I)^{\mathfrak N}=\mathfrak W(R,I)$, 
$\mathfrak D(R,I,I)$ is a nonempty finite set of DVRs, and upon letting
$V_1,\dots,V_h$ be all the distinct members of $\mathfrak D(R,I,I)$, and upon 
letting $W_1,\dots,W_h$ be their respective $(R,I)$-extensions to $L(Z)$, we have 
that $W_1,\dots,W_h$ are distinct DVRs with $E\subset W_j$ for $1\le j\le h$. 

 For $1\le j\le h$ let
$$
P^*_j=E\cap M(W_j)\;\;\text{ and }\;\;\overline P^*_j=E\cap(IW_j)\;\;
\text{ and }\;\;J_j=R\cap(IV_j)
$$ 
but let us postpone considering the rest of the second paragraph.

Turning to the third paragraph: Now
assume that $Y=xZ$ and $A=R[Ix^{-1}]$ where $0\ne x\in R$ is such that
$V_j(x)=V_j(I)$ for all $j$ in a nonempty subset $\Lambda$ of $\{1,\dots,h\}$.
Then for all $j\in\Lambda$ we clearly have $A\subset V_j$ and hence 
$B\subset W_j$. For all $j\in\Lambda$ let
$$
\text{$P_j=A\cap M(V_j)$ with $\overline P_j=A\cap(IV_j)$}
$$
and
$$
\text{$Q_j=B\cap M(W_j)$ with $\overline Q_j=B\cap(IW_j)$.}
$$
Given any $j\in\Lambda$, by taking $V=V_j$ in (6.9) we see that 
$$
\text{$M^*$ is a multiplicative set in $E$ with $E_{M^*}=B$ and
$E\subset B\subset W_j$}
\leqno(1)
$$
and
$$
\begin{cases}
P_j=A\cap Q_j\in\text{spec}(A)\text{ with }Q\in\text{spec}(B)\text{ and }
P^*_j=E\cap Q\in\text{spec}(E)\\
\text{and }P_jB=Q_j=P^*_jB\text{ with }
\text{ht}_AP_j=\text{ht}_BQ_j=\text{ht}_EP^*_j=1.\\
\end{cases}
\leqno(2)
$$
Since $E$ is a normal noetherian domain, by (1) we see that $B$ is also a
normal noetherian domain. Now $B=A[Y,Y^{-1}]$ tells us that
QF$(A)=L\subset L(Y)=\text{QF}(B)$ with $L\cap B=A$, and hence the normality
of $B$ yields the normality of $A$. Thus
$$
\text{all the three rings $A,B,E$ are normal noetherian domains.}
\leqno(3)
$$
By (2) and (3) we see that for all $j\in\Lambda$ we have
$$
\begin{cases}
xA=IA\subset P_j\in\text{nvspec}_A(IA),\\
xB=IB\subset Q_j\in\text{nvspec}_B(IB),\\
P^*_j\in\text{nvspec}_E(IE),\\
A_{P_j}=V_j\text{ and }B_{Q_j}=E_{P^*_j}=W_j.\\
\end{cases}
\leqno(4)
$$
By the last line of (4) we see that for all $j \ne j'$ in $\Lambda$ we have
$$
P_j\ne P_{j'}\;\;\text{ and }\;\;Q_j\ne Q_{j'}\;\;\text{ and }\;\;
P^*_j\ne P^*_{j'}.
\leqno(5)
$$
If $p,\overline p$ is a prime-primary pair (i.e., if $p$ is a prime ideal
and $\overline p$ is a $p$-primary ideal) in a ring, then their contractions
$q,\overline q$ to a subring constitute a prime-primary pair in that subring;
moreover, every nonzero nonunit ideal in a DVR is primary for the maximal ideal;
consequently, for all $j\in\Lambda$ 
$$
\begin{cases}
\text{the ideals $\overline P_j,\overline Q_j,\overline P^*_J$}\\
\text{are $P_j$-primary, $Q_j$-primary, $P^*_j$-primary ideals in $A,B,E$
respectively}
\end{cases}
\leqno(6)
$$
and, because ideals which are primary for distinct prime ideals are obviously
distinct, by (5) we see that for all $j\ne j'$ in $\Lambda$ we have
$$
\overline P_j\ne \overline P_{j'}\;\;\text{ and }\;\;
\overline Q_j\ne \overline Q_{j'}\;\;\text{ and }\;\;
\overline P^*_j\ne \overline P^*_{j'}.
\leqno(7)
$$
In view of (2) to (7) together with (T82) on page 355 of \cite{Ab4}, the 
equation $\mathfrak W(R,I)^{\mathfrak N}=\mathfrak W(R,I)$ implies that 
$$
\begin{cases}
\text{if for all $l\in\{1,\dots,h\}\setminus\Lambda$ we have 
$V_l(x)\ne V_l(I)$}\\
\text{then nvspec$_A(IA)=\{P_j:j\in\Lambda\}$ and}\\
\text{$IA=\cap_{j\in\Lambda}\overline P_j$}\\
\text{is the unique irredundant primary decomposition of $IA$ in $A$.}
\end{cases}
\leqno(8)
$$
In view of (2) to (8) together with (T17) on page 145 and (T30) on page 235
of \cite{Ab4}, the equation $B=A[Y]_{M^*}$ tells us that
$$
\begin{cases}
\text{if for all $l\in\{1,\dots,h\}\setminus\Lambda$ we have 
$V_l(x)\ne V_l(I)$}\\
\text{then nvspec$_B(IB)=\{Q_j:j\in\Lambda\}$ and}\\
\text{$IA=\cap_{j\in\Lambda}\overline Q_j$}\\
\text{is the unique irredundant primary decomposition of $IB$ in $B$.}
\end{cases}
\leqno(9)
$$
Thus we have proved everything in the third paragraph.

\centerline{}

Now let us prove the assertion in the second paragraph which says that
$$
\begin{cases}
\text{$P^*_1,\dots,P^*_h$ are all the distinct members of nvspec$_E(IE)$,}\\ 
\text{ht$_EP^*_j=1$ with $E_{P^*_j}=W_j$ and 
$\overline P^*_j$ is $P^*_j$-primary for $i\le j\le h$,}\\
IE=\overline P^*_1\cap\dots\cap\overline P^*_h\\
\text{is the unique irredundant primary decomposition of $IE$ in $E$}\\
\end{cases}
\leqno(\dagger)
$$
and
$$
I=J_1\cap\dots\cap J_h\text{ with }
J_j=R\cap\overline P^*_j\text{ for }1\le j\le h.
\leqno(\ddagger)
$$
Recall that $R\subset E$ are normal noetherian domains and for $1\le j\le h$ 
we have
$$
\begin{cases}
E\subset W_j\text{ with }P^*_j=E\cap M(W_j)\text{ and }
\overline P^*_j=E\cap(IW_j)\\
\text{and }R\subset V_j\text{ with }J_j=R\cap(IV_j).
\end{cases}
\leqno(1^*)
$$ 
Given any $j$ with $1\le j\le h$, we can clearly find $0\ne x\in I$
with $V_j(x)=V_j(I)$, and then upon taking $Y=xZ$ and $A=R[Ix^{-1}]$
with $\Lambda=\{j\}$, by (2), (3), (4), (6) we see that
$$
\begin{cases}
P^*_j\in\text{nvspec}_E(IE)\text{ with ht}_EP^*_j=1\text{ and }E_{P^*_j}=W_j\\
\text{and $\overline P^*_j$ is a $P^*_j$-primary ideal in $E$.}
\end{cases}
\leqno(2^*)
$$ 

In $(\ddagger)$ it is clear that $I\subset J_1\cap\dots\cap J_h$ with
$J_j=R\cap\overline P^*_j$ for $1\le j\le h$. We shall show that 
$J_1\cap\dots\cap J_h\subset I$ and this will complete the proof of $(\ddagger)$.
Since $I$ is a complete ideal in the normal noetherian domain $R$, 
by definition we have
$$
I=\cap_{V\in\overline D(L/R)}(R\cap I(V))
$$ 
where, for each $V\in\overline D(L/R)$, $I(V)$ is some ideal in $V$; recall 
that $\overline D(L/R)=$ the set of all valuation rings $V$ with QF$(V)=L$ and
$R\subset V$. It follows that
$$
I=\cap_{V\in\overline D(L/R)}(R\cap(IV)).
$$ 
Given any $V\in\overline D(L/R)$ with $I\subset M(V)$ we shall show that 
$J_1\cap\dots\cap J_h\subset R\cap(IV)$ and this will complete the proof of 
$(\ddagger)$. We may assume that $Y=xZ$ and $A=R[Ix^{-1}]$ where $0\ne x\in R$
is such that $V(x)=V(I)$. Then $A\subset V$ and $xA=IA\subset A\cap(IV)$.
Since $\mathfrak W(R,I)^{\mathfrak N}=\mathfrak W(R,I)$, it follows that
$V_j(x)=V_j(I)$. So we may assume that $\Lambda=\{j:1\le j\le h\text{ and }
V_j(x)=V_j(I)\}$. By (8) we have $IA=\cap_{j\in\Lambda}\overline P_j$
and hence 
$$
\cap_{j\in\Lambda}\overline P_j\subset A\cap(IV).
$$
By definition $\overline P_j=A\cap(IV_j)$ and clearly
$\cap_{1\le j\le h}(A\cap(IV_j))\subset\cap_{j\in\Lambda}(A\cap(IV_j))$; 
therefore by the above display we get
$$
\cap_{1\le j\le h}(A\cap(IV_j))\subset A\cap(IV).
$$
Intersecting both sides with $R$ we conclude that 
$J_1\cap\dots\cap J_h\subset R\cap(IV)$. This completes the proof of 
$(\ddagger)$. Since $I$ is a normal ideal, for every $n\in\mathbb N$,
the ideal $I^{n+1}$ is a normal ideal; consequently by $(\ddagger)$ we see that
for every $n\in \mathbb N$ we have
$$
I^{n+1}=\cap_{1\le j\le h}(R\cap(I^{n+1}V_j)).
\leqno(\ddagger_n)
$$

Given any $f\in\cap_{1\le j\le h}\overline P^*_j$, by using $(\ddagger_n)$
we shall show that $f\in IE$ and, in view of (1*) and (2*), this will prove
$(\dagger)$ which will complete the proof of Lemma (6.11). By (4.1) we can
express $f$ as a finite sum
$f=\sum_{n\in\mathbb N}f_nZ^n$ with $f_n\in I^n$.
By the definition of $W_j$ as the $(R,I)$-extension of $V_j$ we see that
$\overline P^*_j$ is a homogeneous ideal in the homogeneous ring $E$;
alternatively this follows because $\overline P^*_j=E\cap(IE_{P^*_j})=$
the primary component of the homogeneous ideal $IE$ with respect to its
minimal prime $P^*_j$. Therefore $f_nZ^n\in\cap_{1\le j\le h}\overline P^*_j$
for all $n\in\mathbb N$, and it suffices to show that for every $n\in\mathbb N$
we have $f_nZ^n\in IE$. Assuming $f_n\ne 0$, for $1\le j\le h$ we clearly have
$$
f_nZ^n\in\overline P^*_j\Rightarrow W_j(f_nZ^n)\ge V_j(I)
$$
with
$$
W_j(f_nZ^n)=V_j(f_n)-nV_j(I)
$$
and hence
$$
f_nZ^n\in\overline P^*_j\Rightarrow V_j(f_n)-nV_j(I)\ge V_j(I)
\Rightarrow V_j(f_n)\ge (n+1)V_j(I)=V_j(I^{n+1}).
$$
Therefore
$$
f_nZ^n\in\cap_{1\le j\le h}\overline P^*_j\Rightarrow 
f_n\in\cap_{1\le j\le h}(R\cap(I^{n+1}V_j))
$$
and hence by $(\ddagger_n)$ we conclude that
$$
f_nZ^n\in\cap_{1\le j\le h}\overline P^*_j\Rightarrow 
f_n\in I^{n+1}
$$
and clearly
$$
f_n\in I^{n+1}\Rightarrow f_nZ^n\in IE.
$$

\centerline{}

{\it Remark $(6.11^\flat)$.} 
Alternatively, the normality of $A$ and $B$ can be seen thus.  $A$ is normal 
because $\mathfrak V(A)$ is an affine piece of the normal variety 
$\mathfrak W(R,I)$. Therefore $B$ is normal because it is a localization of
$A[Y]$.

\centerline{}

DEFINITION-OBSERVATION (6.12). For any nonunit ideal $J$ in a ring $S$ and
any nonnegative integer $i$ we define the {\bf depth $i$ portion} of
nvspec$_SJ$ by putting
$$
(\text{nvspec}_SJ)_i=\{P\in\text{nvspec}_SJ:\text{dpt}_SJ=i\}.
$$

To use the above definition, referring to pages 206-215 and 399-408 of 
\cite{Ab4} for the details of the theory of homogeneous rings and irrelevant 
ideals, let 
$$
F=\sum_{n\in\mathbb N}F_n
$$ 
be a homogeneous ring over a field $F_0$ with $[F_1:F_0]<\infty$.
As usual let
$$
\text{$\Omega(F)=F_1F=\sum_{n\in\mathbb N_+}F_n=$ 
the unique homogeneous maximal ideal in $F$}
$$
and let
$$
\overline\Omega(F)=\cup_{n\in\mathbb N_+}F_n.
$$
Recall that a nonunit homogeneous ideal $G$ in $F$ is irrelevant means 
$\Omega(F)\subset\text{rad}_FG$, i.e., equivalently, $\Omega(F)=\text{rad}_FG$. 
Note that dim$(F)\in\mathbb N$ and for any nonunit homogeneous ideal $G$ in
$F$ we have dim$(F/G)\in\mathbb N$ with $\dim(F/G)\le\dim(F)$.

Now given any nonunit homogeneous ideal $G$ in $F$, upon letting
$\dim(F/G)=e$, let us prove the following Observations (I) to (VI).

(I) $F$ is integral over $F_0[G]\Leftrightarrow G$ is irrelevant 
$\Leftrightarrow e=0$.

(II) $\emptyset\ne(\text{nvspec}_FG)_e\subset\text{nvspec}_FG$.
If $e\le 1$ then $(\text{nvspec}_FG)_e=\text{nvspec}_FG$.

(III) For any $y\in\overline\Omega(F)$ we have: 
$\dim(F/(G\cup\{y\})F)=e-1\Leftrightarrow y\not\in P$
for all $P\in(\text{nvspec}_FG)_e$.

(IV) For any $y\in\overline\Omega(F)$ we have dim$(F/(G\cup\{y\})F)\ge e-1$.

(V) For any elements $y_1,\dots,y_e$ in $\overline\Omega(F)$,
upon letting $G_i=(G\cup\{y_1,\dots,y_i\})F$ for $0\le i\le e$,
the following conditions (1) to (4) are mutually equivalent.

(1) For $1\le i\le e$ we have $\dim(F/G_i)=e-i$.

(2) For $1\le i\le e$ we have $\dim(F/G_i)\le e-i$.

(3) For $1\le i\le e$ we have 
$y_i\not\in P$ for all $P\in(\text{nvspec}_FG_{i-1})_{e-i+1}$.

(4) For $1\le i\le e$ we have 
$y_i\not\in P$ for all $P\in \widehat G_{i-1}$ where
$\widehat G_{i-1}=\text{nvspec}_FG_{i-1}$ or
$\widehat G_{i-1}=(\text{nvspec}_FG_{i-1})_{e-i+1}$ according as
$i=e$ or $i\ne e$.

(VI) There exist elements $z_1,\dots,z_e$ in $F_s$ for some $s\in\mathbb N_+$
such that upon letting 
$$
G_e=(G\cup\{z_1,\dots,z_e\})F
$$
we have that dim$(F/G_i)=0$ and $F$ is integral over $F_0[G_e]$.
Moreover, if $F_0$ is infinite then we can take $s=1$. 

\centerline{}

PROOF. Let us observe that the associated primes of any homogeneous ideal are 
homogeneous, and $F_1F$ is the only maximal ideal in $F$ which is homogeneous.

Now (I) follows from (T104) on page 401 of \cite{Ab4}. 

The proofs of (II) and (III) are straightforward.

To prove (IV), consider the homogeneous ring $F^*=F/G$. Also consider the
homogeneous ring $F'=F/(G\cup\{y\})F$ and let $\dim(F')=e'$.  By the 
homogeneous noether normalization theorem (T106) on page 408 of \cite{Ab4}, 
we can find elements $z_1,\dots,z_{e'}$ in $\overline\Omega(F)$ whose images 
$z'_1,\dots,z'_{e'}$ in $F'$ are such that $F'$ is integral over
$F'_0[z'_1,\dots,z'_{e'}]$. By (I) we see that $(z'_1,\dots,z'_{e'})F'$
is an irrelevant ideal in $F'$, and hence upon letting
$y^*,z^*_1,\dots,z^*_{e'}$ be the respective images of $y,z_1,\dots,z_{e'}$
in $F^*$ it follows that $(y^*,z^*_1,\dots,z^*_{e'})F^*$ is an irrelevant ideal
in $F^*$. Therefore, again by (I), $F^*$ is integral over
$F^*_0[y^*,z^*_1,\dots,z^*_{e'}]$. Consequently, say by
(O10), (O11), (T45), (T47) on pages 110, 111, 247, 250 of \cite{Ab4},
we get $1+e'\ge e$ and hence $e'\ge e-1$.

Turning to (V), by (IV) we get $(1)\Leftrightarrow(2)$,
by (III) we get $(1)\Leftrightarrow(3)$,
and by (II) and (III) we get $(1)\Leftrightarrow(4)$.

To prove (VI), again consider the homogeneous ring $F^*=F/G$.
Now by the normalization theorem (T46) and the homogeneous 
normalization theorem (T106) respectively on pages 248 and 408 of \cite{Ab4},
we can find elements $z_1,\dots,z_e$ in $F_s$ for some $s\in\mathbb N_+$, 
where we can take $s=1$ in case $F_0$ is infinite, such that upon letting
$z_1^*,\dots,z^*_e$ be their respective images in $F^*$ we have that
$F^*$ is integral over $F^*_0[z^*_1,\dots,z^*_e]$. By (I) it follows that 
$(z^*_1,\dots,z^*_e)F^*$ is an irrelevant ideal in $F^*$. Therefore
$$
G_e=(G\cup\{z_1,\dots,z_e\})F
$$ 
is an irrelevant ideal in $F$. Therefore by (I) we conclude that 
dim$(F/G_e)=0$ and $F$ is integral over $F_0[G_e]$.
 
\centerline{}

LEMMA (6.13). Let the assumptions be as in (6.11). Also assume that $R$ is a 
$d$-dimensional normal local domain with $M=M(R)$,
and $I$ is a nonzero normal $M$-primary ideal in $R$. Let 
$$
\text{$F=F_R(I)=E/ME=\sum_{n\in\mathbb N}F_n=$ the form ring of $I$}
$$
and let
$$
\mu_n:I^n\to F_n\;\;\text{ with }\;\; \text{ker}(\mu_n)=MI^n
$$
be the canoninal $R$-epimorphism. For $1\le j\le h$ let
$$
Q'_j=\mu_1(P^*_j)\;\;\text{ and }\;\; P'_j=Q'_jF.
$$
Also let
$$
\Lambda'=\{j:1\le j\le h\text{ and dpt}_EP^*_j=d\}.
$$
Then we have the following,

(I) $P^*_1,\dots,P^*_h$ are all the distinct members of nvspec$_E(ME)$
and we have
$$
\text{ht}_EP^*_1=\dots=\text{ht}_EP^*_h=1
$$
and
$$
\text{max}(\text{dpt}_EP^*_1,\dots,\text{dpt}_EP^*_h)=d=\dim(E)-1.
$$ 

(II) $P'_1,\dots,P'_h$ are all the distinct members of nvspec$_F\{0\}$ and we have
$$
\text{ht}_FP'_1=\dots=\text{ht}_FP'_h=0
$$
and
$$
\text{max}(\text{dpt}_FP'_1,\dots,\text{dpt}_FP'_h)=d=\dim(F).
$$ 
Also we have
$$
\Lambda'=\{j:1\le j\le h\text{ and dpt}_FP'_j=d\}.
$$

(III) Given any $j\in\{1,\dots,h\}$ and any $x\in I$ we have
$$
V_j(x)=V_j(I)\Leftrightarrow xZ\not\in P^*_j\Leftrightarrow\mu_1(x)\not\in P'_j.
$$

(IV) Given any elements $x_1,\dots,x_d$ in $I$ let 
$y_1=\mu_1(x_1),\dots,y_d=\mu_1(x_d)$. Then $y_1,\dots,y_d$ are elements in $F_1$
and we have
$$
\text{$(x_1,\dots,x_d)R$ is a reduction of 
$I\Leftrightarrow F$ is integral over $F_0[y_1,\dots,y_h]$.}
$$

(V) If $x_1,\dots,x_d$ are elements in $I$ such that $(x_1,\dots,x_d)R$ is a
reduction of $I$ then for $1\le i\le d$ we have $V_j(x_i)=V_j(I)$
for all $j\in\Lambda'$.

(VI) If $d=1$ then $h=1$ and $\Lambda'=\{1\}$.

(VII) If $d=2$ then $\Lambda'=\{1,\dots,h\}$.

(VIII) If $d=1$ then for any $x_1\in I$ we have that
$$
V_1(x)=V_1(I)\Leftrightarrow x_1R\text{ is a reduction of }I.
$$

(IX) If $R/M$ is infinite 
then there exist elements $x_1,\dots,x_d$ in $I$ 
such that $(x_1,\dots,x_d)R$ is a reduction of $I$.

(X) If $R/M$ is infinite and $x_1\in I$ is such that 
$V_j(x_1)=V_j(I)$ for all $j\in\Lambda'$ then there exist elements 
$x_2,\dots,x_d$ in $I$ such that $(x_1,\dots,x_d)R$ is a reduction of $I$.

(XI) Given any $j\in\{1,\dots,h\}$ there exists $x'_j\in I$ such that 
$V_j(x'_j)=V_j(I)$.  If $R/M$ is infinite 
then there exists $x_1\in I$ such that $V_j(x_1)=V_j(I)$ for $1\le j\le h$.

\centerline{}

PROOF. By (6.11) we know that
$$
\text{ht}_EP^*_1=\dots=\text{ht}_EP^*_h=1
$$
and $P^*_1,\dots,P^*_h$ are all the distinct members of nvspec$_E(IE)$.
Since $I$ is $M$-primary, it follows that
$P^*_1,\dots,P^*_h$ are all the distinct members of nvspec$_E(ME)$.
By (6.4) of \cite{AbH} we also have $\dim(F)=d$ and hence we get (II).
Now it follows that
$$
\text{max}(\text{dpt}_EP^*_1,\dots,\text{dpt}_EP^*_h)=d\ge \dim(E)-1.
$$ 
So to complete the proof of (I) and (II) we only need to show that 
dim$(E)\le \dim(R)+1$. But this follows from the Multiple Ring Extension 
Lemma (T55) on page 269 of \cite{Ab4} by noting that $E$ is an affine domain
over $R$ and transcendence degree of QF$(E)$ over the quotient field of
$R$ is $1$.

The second implication in (III) follows from the fact that
$\mu_1(x)$ is the image of $xZ$ under the residue class epimorphism $E\to F$.
The first implication of (III) follows by noting that for $0\ne x\in I$
we clearly have 
$$
V_j(x)=V_j(I)\Leftrightarrow V_j(x)\le V_j(I)
$$
and 
$$
xZ\not\in P^*_j\Leftrightarrow W_j(xZ)\le 0
$$
and
$$
W_j(xZ)=V_j(x)-W_j(Z)=V_j(x)-V_j(I).
$$

(IV) follows from (6.1) of \cite{AbH}.

In view of (6.12), assertions (V) to (X) follow from assertions (I) to (IV)
where we note that: in proving (IX) we take $G=\{0\}F$ in (6.12)(VI),
while in proving (X) we take $G=\mu_1(x_1)F$ in (6.12)(VI).

The first part of (XI) is obvious, and from it to deduce the second part,
assume that the field $R/M$ is infinite. Let $x_1=a_1x'_1+\dots+a_hx'_h$
where $a_1,\dots,a_h$ in $R$ are to be chosen. For any $a\in R$ let 
$\overline a$ be its image in $R/M$. Clearly 
$V_j(x_1)>V_j(I)\Leftrightarrow(\overline a_1,\dots,\overline a_h)$ belongs to 
a certain proper subspace $K_j$ of $(R/M)^h$. 
The infiniteness of $R/M$ implies that $K_1\cup\dots\cup K_h\ne(R/M)^h$ and it
suffices to take $a_1,\dots,a_h$ to be such that
$(\overline a_1,\dots,\overline a_h)\not\in K_1\cup\dots\cup K_h$.

\centerline{}

Without assuming $R/M$ to be infinite, we shall now prove the following
variation of parts (IX) to (XI) of (6.13).

\centerline{}

LEMMA (6.14). Let the assumptions be as in (6.13). Then, without assuming $R/M$ 
to be infinite, we have the following.

(I) Given any $r\in\mathbb N_+$, there exist elements $x_1,\dots,x_d$ in 
$I^{rs}$ for some $s\in\mathbb N_+$ 
such that, for every $t\in\mathbb N_+$,
$(x^t_1,\dots,x^t_d)R$ is a reduction of $I^{rst}$.

(II) If $r\in\mathbb N_+$ and $x_1\in I^r$ are such that 
$V_j(x_1)=V_j(I^r)$ for all $j\in\Lambda'$ then there exist elements 
$x_2,\dots,x_d$ in $I^{rs}$ for some $s\in\mathbb N_+$ such that, for every 
$t\in\mathbb N_+$, $(x^{st}_1,x_2^t,\dots,x^t_d)R$ is a reduction of $I^{rst}$.

(III) Given any $r\in\mathbb N_+$, there exist $s\in\mathbb N_+$ and 
$x_1\in I^{rs}$ such that 
$V_j(x_1)=V_j(I^{rs})$ for all $j\in\Lambda'$.

\centerline{}

PROOF. Given any $q\in\mathbb N_+$, clearly $I^q$ is a nonzero normal
$M$-primary ideal in $R$ with $\mathfrak W(R,I^q)=\mathfrak W(R,I)$. It
follows that $V_1,\dots,V_h$ are all the distinct members of
$\mathfrak D(R,I^q,I^q)$. Moreover 
$$
E^{(q)}=E_R(I^q)\subset E_R(I)=E
$$
and $W_1,\dots,W_h$ are the respective $(R,I^q)$-extensions of $V_1,\dots,V_h$.
Furthermore 
$$
P^{(q)}_j=E^{(q)}\cap P^*_j=E^{(q)}\cap M(W_j)\in\text{spec}(E^{(q)})
\;\;\text{ with }\;\;\text{ht}_{E^{(q)}}P^{(q)}_j=1
$$
for $1\le j\le h$, and the prime ideals $P^{(q)}_1,\dots,P^{(q)}_h$ 
are all the distinct members of nvspec$_{E^{(q)}}(I^qE^{(q)})$. Finally 
$\Lambda'=\{j:1\le j\le h\text{ and }\text{dpt}_{E^{(q)}}P^{(q)}_j=d\}$. 
Note that the above two displays 
include the definitions of the symbols $E^{(q)}$ and $P^{(q)}_j$. Let
$$
\text{$F^{(q)}=F_R(I^q)=E^{(q)}/ME^{(q)}=\sum_{n\in\mathbb N}F^{(q)}_n=$ 
the form ring of $I^q$}
$$
and let
$$
\mu^{(q)}_n:I^n\to F^{(q)}_n\;\;\text{ with }\;\; \text{ker}(\mu^{(q)}_n)=MI^n
$$
be the canonical $R$-epimorphism. Applying (4.1) to $E$ and $E^{(q)}$ we see
that
$$ 
ME^{(q)}=E^{(q)}\cap(ME)
$$
and hence, upon letting
$$ 
\psi^{(q)}:E^{(q)}\to E
$$
be the inclusion monomorphism and
$$
\nu:E\to F\;\;\text{ and }\;\;\nu^{(q)}:E^{(q)}\to F^{(q)}
$$
be the residue class epimorphisms, there exists a unique monomorphism
$$
\phi^{(q)}:F^{(q)}\to F
$$
such that
$$
\phi^{(q)}\nu^{(q)}=\nu\psi^{(q)}.
$$
Note that the restriction of $\psi^{(q)}$ to $E^{(q)}_n$ gives an isomorphism
$E^{(q)}_n\to E_{nq}$, and the restriction of $\phi^{(q)}$ to $F^{(q)}_n$ gives 
an isomorphism $F^{(q)}_n\to F_{nq}$; these isomorphisms are the foundation of
the {\bf Veronese Embedding} expounded on pages 263-283 of \cite{Ab3};
let it be recorded that VERONESE was the Param-Param-Guru of Abhyankar whose
Guru Zariski was a pupil of Castelnuovo whose Guru was Veronese; going back one 
step further, Cremona was the Guru of Veronese; this makes Cremona Abhyankar's
{\bf Param-Param-Param-Guru}. The commutative diagram
\begin{equation*}
\CD
E_{nq} @>>> E @>{\nu}>> F @<<< F_{nq}  \\
@AAA  @A{\psi^{(q)}}AA    @A{\phi^{(q)}}AA  @AAA \\
E_n^{(q)} @>>> E^{(q)}  @>{\nu^{(q)}}>> F^{(q)} @<<<  F_n^{(q)}.
\endCD
\end{equation*}
exhibits the various maps which we have discussed.

The above observations will be used tacitly. 

In proving (I) to (III) we shall assume that $r=1$; the general case will
then follow by taking $I^r$ for $I$.

To prove (I) and (II), by (6.12)(VI) we can find elements $z_1,\dots,z_d$
in $F_s$ for some $s\in\mathbb N_+$ such that the ring $F$ is integral over the
subring $F_0[z_1,\dots,z_d]$, where in case of (I) we take
$G=\{0\}F$, while in case of (II) we take
$G=\mu_1(x_1)$ and $z_1=\mu_1(x_1)^s$.
Using the monomorphism $\phi^{(s)}$ we get unique elements $y_1,\dots,y_d$
in $F^{(s)}_1$ such that 
$\phi^{(s)}(y_1)=z_1,\dots,\phi^{(s)}(y_d)=z_d$, and then the ring $F^{(s)}$ 
is clearly integral over the subring $F^{(s)}_0[y_1,\dots,y_d]$.
In case of (I) we can take elements $x_1,\dots,x_d$ in $I^s$ such that
$\mu^{(s)}_1(x_1)=y_1,\dots,\mu^{(s)}_1(x_d)=y_d$, while in case of (II) we can 
take elements $x_2,\dots,x_d$ in $I^s$ such that 
$\mu^{(s)}_1(x_2)=y_2,\dots,\mu^{(s)}_1(x_d)=y_d$. By (6.13) we see that, in case 
of (I), $(x_1,\dots,x_d)R$ is a reduction of $I^s$, while, in case of (II),
$(x_1^s,x_2,\dots,x_d)R$ is a reduction of $I^s$. Given any $t\in\mathbb N_+$,
using (6.13)(V) and the monomorphisms $\phi^{(s)}$ and $\phi^{(st)}$,
we conclude that, in case of (I), $(x_1^t,\dots,x_d^t)R$ is a reduction of 
$I^{st}$, while, in case of (II), 
$(x_1^{st},x_2^t,\dots,x_d^t)R$ is a reduction of $I^{st}$. 

This proves (I) and (II). In view of (6.13)(V), (III) follows from (I).

\centerline{}

\centerline{}

{{\bf Section 7: Extended Rees Rings.}} 
An alternative way of approaching parts of (6.11) to
(6.14) is provided by the theory of extended Rees rings. To introduce 
these rings, let $I$ be an ideal in a nonnull ring $R$.
The {\bf extended Rees ring} $\widehat E_R(I)$ of $I$ relative to $R$ with
variable $Z$ is defined by putting
$$
\widehat E_R(I)=R[Z^{-1},IZ]=\sum_{n\in\mathbb Z}\widehat E_R(I)_n
$$
which makes it a $\mathbb Z$-graded = an integrally graded ring.
Note that $E_R(I)$ is a graded subring of $\widehat E_R(I)$ and we have
$$
\widehat E_R(I)_n=\begin{cases}
\{gZ^n:g\in R\}&\text{if }n<0\\
\{gZ^n:g\in I^n\}=E_R(I)_n&\text{if }n\ge 0.\end{cases}
$$
At any rate, every $f\in\widehat E_R(I)$ can uniquely be written as a finite sum
$$
f=\sum_{n\in\mathbb Z}f_nZ^n\;\;\text{ with }\;\;
f_n\in\begin{cases}R&\text{ if }n<0\\I^n&\text{ if }n\ge 0.\end{cases}
\leqno(7.1)
$$
Now consider
$$
E=E_R(I)=\sum_{n\in\mathbb N}E_n\subset
\widehat E=\widehat E_R(I)=\sum_{n\in\mathbb Z}\widehat E_n
\subset R'=R[Z^{-1},Z]=\sum_{n\in\mathbb Z}R'_n
$$
and note that
$$
R'=R^-\oplus R^+\;\;\text{ where }\;\;
R^-=\sum_{n\in\mathbb Z\setminus\mathbb N}R^-_n
\;\;\text{ and }\;\;R^+=\sum_{n\in\mathbb N}R^+_n.
$$
By (4.1) and (7.1) we see that
$$
\begin{cases}
E\cap(Z^{-1}\widehat E)=IE\\
\text{which induces an isomorphism}\\
\widehat E/(Z^{-1}\widehat E)\approx E/(IE).
\end{cases}
\leqno(7.2)
$$
Heuristically speaking, (7.2) says that $Z^{-1}$ in $\widehat E$ is sort of a
generic element of $I$, but $I\widehat E$ is properly contained in $Z^{-1}E$.
Indeed by (4.1) and (7,1) we see that
$$
R^-\cap(Z^{-1}\widehat E)=R^-\;\;\text{ and }\;\;
R^+\cap(Z^{-1}\widehat E)=IE
\leqno(7.3)
$$
i.e., the negative portion of $Z^{-1}\widehat E$ equals the entire $R^-$, and 
the nonnegative portion of $Z^{-1}\widehat E$ equals $IE$.
By (4.1) and (7.1) we also see that
$$
\begin{cases}
\text{for any homogeneous ideal $J$ in $E$,
upon letting $\widehat J=R^-\oplus J$,}\\
\text{$\widehat J$ is a homogeneous ideal in $\widehat E$ such that}\\
R^-\cap\widehat J=R^-\;\;\text{ and }\;\;R^+\cap\widehat J=J.\\
\end{cases}
\leqno(7.4)
$$
Taking $E=J$ in (7.4) we get
$$
\widehat E=R^-\oplus E\;\;\text{ with }\;\;R^-\cap\widehat E=R^-
\;\;\text{ and }\;\;R^+\cap\widehat E=E.
\leqno(7.5)
$$
Let us now prove a lemma about $\widehat E$.

\centerline{}

LEMMA (7.6). Let the assumptions be as in (6.11).
Then $\widehat E$ is a normal noetherian domain and, upon letting 
$P'_1,\dots,P'_{h'}$ be the minimal primes of $Z^{-1}\widehat E$ and upon 
letting $W'_j=\widehat E_{P'_j}$ for $1\le j\le h'$, we have the following.
[Note that, by Krull Normality Lemma (T82) on page 355 of \cite{Ab4}, 
all associated primes of $Z^{-1}\widehat E$ are minimal and have height one].

(I) Upon letting $M'$ be the multiplicative set $\{1,Z^{-1},Z^{-2},\dots\}$
in $\widehat E$ we have
$$
\text{$\widehat E_{M'}=R[Z^{-1},Z]$
\;\;and hence\;\; $\widehat E=R[Z^{-1},Z]\cap W'_1\cap\dots\cap W'_{h'}$.}
$$

(II) $h=h'$ and, after a suitable (obviously unique)
relabelling, for $1\le j\le h$ we have 
$$
P'_j=\widehat{P^*_j}\;\;\text{ and }\;\; W'_j=W_j.
$$

\centerline{}

PROOF.
The Krull Normality Lemma (T82) on page 355 of \cite{Ab4} gives us (I). 
Since $P'_1,\dots,P'_{h'}$ are the minimal primes of $Z^{-1}\widehat E$ in
$\widehat E$ and $P^*_1,\dots,P^*_h$ are the minimal primes of $IE$ in $E$,
by (7.2) and (7.4) we see that $h'=h$ and after a suitable relabelling we have
$P'_j=\widehat{P^*_j}$ for $1\le j\le h$. Now, since $\widehat E_{P'_j}=W'_j$
and $\widehat E\subset E_{P^*_j}=W_j$, we get $W'_j=W_j$ for $1\le j\le h$.

\centerline{}

{\it Remark $(7.6^\flat)$.} The purport of Lemma (7.6) is that we can either
first get the $P'_j$ and the $W'_j$ and then obtain $W_j$ and $P^*_j$ satisfying
(6.11), or we can first do (6.11) and then get hold of the $P'_j$ and the
$W'_j$ belonging to $\widehat E$ by tacking on the negative piece $R^-$. The
ideal $Z^{-1}\widehat E$ contains $R^-$ and is hence quite powerful.

\centerline{}

{{\bf Section 8: Dicriticals of Two Dimensional Regular Local Domains.}} 
We are now ready to reap the harvest 
from the work done in the previous Sections. 

Let $R$ be a $d$-dimensional local domain with quotient field $L$, let
$I$ be a nonzero ideal in $R$ with $I\subset M=M(R)$, let $Z$ be an 
indeterminate over the quotient field $L$ of $R$, and let 
$$
\text{$E=E_R(I)=R[IZ]=$ the Rees ring of $I$ relative to $R$ with variable 
$Z$.}
$$ 
Let
$$
\text{$\mathfrak D(R,I)=(\mathfrak W(R,I)^{\Delta}_1)^{\mathfrak N}=$
the set of all dicritical divisors of $I$ in $R$.}
$$
 From (6.12)(III) recall that then: 
$\mathfrak D(R,I)$ is a finite set of DVRs; upon letting 
$V_1,\dots,V_h$ to be all the distinct members of $\mathfrak D(R,J)$, 
and upon letting $W_1,\dots,W_h$ be their respective
$(R,I)$-extensions to $L(Z)$, we have that $W_1,\dots,W_h$ are
distinct DVRs with $E\subset W_j$ for $1\le j\le h$; for the definition of
$(R,I)$-extension see the second paragraph of Section 6;
note that if $I$ is $M$-primary then $h>0$.
Let
$$
P^*_j=E\cap M(W_j)\in\text{spec}(E)\;\;\text{ for }\;\;1\le j\le h
$$
and let
$$
\Lambda'=\{j:1\le j\le h\text{ and dpt}_EP^*_j=d\}.
$$

\centerline{}

THEOREM (8.1). Assuming that $R$ is a $d$-dimensional normal local domain and 
$I$ is a nonzero normal $M$-primary ideal in $R$, we have the following.

(I) $P^*_1,\dots,P^*_h$ are all the distinct members of nvspec$_E(ME)$
and we have
$$
\text{ht}_EP^*_1=\dots=\text{ht}_EP^*_h=1
$$
and
$$
\text{max}(\text{dpt}_EP^*_1,\dots,\text{dpt}_EP^*_h)=d=\dim(E)-1.
$$ 

(II) Given any $j\in\{1,\dots,h\}$ and any $x\in I$ we have
$$
V_j(x)=V_j(I)\Leftrightarrow xZ\not\in P^*_j.
$$

(III) If $x_1,\dots,x_d$ are elements in $I$ such that $(x_1,\dots,x_d)R$ is a
reduction of $I$ then for $1\le i\le d$ we have $V_j(x_i)=V_j(I)$
for all $j\in\Lambda'$.

(IV) If $d=1$ then $h=1$ and $\Lambda'=\{1\}$.
If $d=2$ then $\Lambda'=\{1,\dots,h\}$.

(V) If $d=1$ then for any $x_1\in I$ we have that
$$
V_1(x)=V_1(I)\Leftrightarrow x_1R\text{ is a reduction of }I.
$$

(VI) Given any $r\in\mathbb N_+$, there exist elements $x_1,\dots,x_d$ in 
$I^{rs}$ for some $s\in\mathbb N_+$ 
such that, for every $t\in\mathbb N_+$,
$(x^t_1,\dots,x^t_d)R$ is a reduction of $I^{rst}$.

(VII) If $r\in\mathbb N_+$ and $x_1\in I^r$ are such that 
$V_j(x_1)=V_j(I^r)$ for all $j\in\Lambda'$ then there exist elements 
$x_2,\dots,x_d$ in $I^{rs}$ for some $s\in\mathbb N_+$ such that, for every 
$t\in\mathbb N_+$, $(x^{st}_1,x_2^t,\dots,x^t_d)R$ is a reduction of $I^{rst}$.

(VIII) Given any $r\in\mathbb N_+$, there exist $s\in\mathbb N_+$ and 
$x_1\in I^{rs}$ such that 
$V_j(x_1)=V_j(I^{rs})$ for all $j\in\Lambda'$.

(IX) If $R/M$ is infinite 
then there exist elements $x_1,\dots,x_d$ in $I$ 
such that $(x_1,\dots,x_d)R$ is a reduction of $I$.

(X) If $R/M$ is infinite and $x_1\in I$ is such that 
$V_j(x_1)=V_j(I)$ for all $j\in\Lambda'$ then there exist elements 
$x_2,\dots,x_d$ in $I$ such that $(x_1,\dots,x_d)R$ is a reduction of $I$.

(XI) Given any $j\in\{1,\dots,h\}$ there exists $x'_j\in I$ such that 
$V_j(x'_j)=V_j(I)$.  If $R/M$ is infinite 
then there exists $x_1\in I$ such that $V_j(x_1)=V_j(I)$ for $1\le j\le h$.

\centerline{}

PROOF. We are done by (6.12) to (6.14).

\centerline{}

{\it Remark $(8.1^\flat)$.} The form ring $F_R(I)$ was used as a tool in
proving Theorems (8.1) and (8.2) but does not explicitly appear in their 
statements. Similarly, the Rees ring $E_R(I)$ was used as a tool 
in proving these theorems and does appear in their statements, but it is
not referred to in parts (V), (VI), (IX), (X), (XI) of Theorem (8.1) and in
parts (V) to (XI) of Theorem (8.2); in this list we can include part (VII) 
of Theorem (8.1) by changing the phrase ``for all $j\in\Lambda'$ then''
by the phrase ``for $1\le j\le h$ then''.

\centerline{}

THEOREM (8.2). Either assume that $R$ is a two dimensional regular local 
domain and $I$ is a complete $M$-primary ideal in $R$, or assume that $R$ is
a two dimensional normal local domain and $I$ is a normal $M$-primary ideal
in $R$. Then we have the following.

(I) $P^*_1,\dots,P^*_h$ are all the distinct members of nvspec$_E(ME)$
and we have
$$
\text{ht}_EP^*_1=\dots=\text{ht}_EP^*_h=1.
$$

(II) Given any $j\in\{1,\dots,h\}$ and any $x\in I$ we have
$$
V_j(x)=V_j(I)\Leftrightarrow xZ\not\in P^*_j.
$$

(III) If $x_1,x_2$ are elements in $I$ such that $(x_1,x_2)R$ is a
reduction of $I$ then for $1\le i\le 2$ and $1\le j\le h$ we have 
$V_j(x_i)=V_j(I)$.

(IV) We have
$$
\text{dpt}_EP^*_1=\dots=\text{dpt}_EP^*_h=2.
$$ 

(V) We have dim$(E)=3$.

(VI) Given any $r\in\mathbb N_+$, there exist elements $x_1,x_2$ in 
$I^{rs}$ for some $s\in\mathbb N_+$ 
such that, for every $t\in\mathbb N_+$,
$(x^t_1,x^t_2)R$ is a reduction of $I^{rst}$.

(VII) If $r\in\mathbb N_+$ and $x_1\in I^r$ are such that 
$V_j(x_1)=V_j(I^r)$ for $1\le j\le h$ then there exists 
$x_2\in I^{rs}$ for some $s\in\mathbb N_+$ such that, for every 
$t\in\mathbb N_+$, $(x^{st}_1,x^t_2)R$ is a reduction of $I^{rst}$.

(VIII) Given any $r\in\mathbb N_+$, there exist $s\in\mathbb N_+$ and 
$x_1\in I^{rs}$ such that 
$V_j(x_1)=V_j(I^{rs})$ for $1\le j\le h$.

(IX) If $R/M$ is infinite 
then there exist elements $x_1,x_2$ in $I$ 
such that $(x_1,x_2)R$ is a reduction of $I$.

(X) If $R/M$ is infinite and $x_1\in I$ is such that 
$V_j(x_1)=V_j(I)$ for $1\le j\le h$ then there exists 
$x_2\in I$ such that $(x_1,x_2)R$ is a reduction of $I$.

(XI) Given any $j\in\{1,\dots,h\}$ there exists $x'_j\in I$ such that 
$V_j(x'_j)=V_j(I)$.  If $R/M$ is infinite 
then there exists $x_1\in I$ such that $V_j(x_1)=V_j(I)$ for $1\le j\le h$.

\centerline{}

PROOF.  Zariski's Theorem $(2')$ on page 385 of Volume II of \cite{Zar} tells
us that if $R$ is two dimensional regular local domain then every complete
$M$-primary ideal in $R$ is a normal ideal in $R$.
Consequently we are done by Theorem (8.1).

\centerline{}

{\it Remark $(8.2^\flat)$.} Proposition (3.5) of \cite{Ab11} says the following.

\centerline{}

$(\dagger)$ Assume that $R$ is a two dimensional regular local domain.
For any $V,W$ in $D(R)^\Delta$ let $c(R,V,W)=\ord_V\zeta_R(W)$, 
i.e., $c(R,V,W)=\text{min}\{\ord_V\theta:\theta\in\zeta_R(W)\}$.
Let $J$ be a special primary pencil at $R$ and assume that $J=(a,b)R$ where
$b=\eta^m$ with $\eta\in M(R)\setminus M(R)^2$ and $m\in\mathbb N_+$.
Let $\mathfrak D(R,J)=U$. Then clearly there exists
$n(W)\in\mathbb N_+$ for all $W\in U$ such that 
$$
I=\prod_{W\in U}\zeta_R(W)^{n(W)}
\leqno(\bullet)
$$
is the integral closure of $J$ in $R$, and hence
$$
m~\ord_V\eta=\sum_{W\in U}n(W)c(R,V,W)
\;\;\text{ for all $V\in U$.}
\leqno(\bullet\bullet)
$$

\centerline{}

As said in the Introductory Note (3.5)(0) of \cite{Ab11}, the proof
of $(\dagger)$ is contained in its statement. Now by (8.2) we get
the converse of $(\dagger)$ stated below.

\centerline{}

$(\ddagger)$ Assume that $R$ is a two dimensional regular local domain.
For any $V,W$ in $D(R)^\Delta$ let $c(R,V,W)=\ord_V\zeta_R(W)$, 
i.e., $c(R,V,W)=\text{min}\{\ord_V\theta:\theta\in\zeta_R(W)\}$. Given any 
nonempty finite $U\subset D(R)^\Delta$, let $I$ be the complete $M$-primary 
ideal in $R$ given by $(\bullet)$ with $n(W)\in\mathbb N_+$. 
Then we have the following.

(1) There exists $b\in I^s$ for some $s\in\mathbb N_+$ such that
for all $V\in U$ we have $V(b)=V(I^s)$; given any such $b$ and $s$, their exists
$a\in I^{st}$ for some $t\in\mathbb N_+$ such that the pencil $J=(a,b^t)R$
is a reduction of $I^{st}$. Moreover, if $R/M$ is infinite then we can take
$s=t=1$.

(2) If $b=\eta^m$ with $\eta\in M(R)\setminus M(R)^2$ and $m\in\mathbb N_+$
are such that $(\bullet\bullet)$ is satisfied, then there exists $a\in I^t$ for 
some $t\in\mathbb N_+$ such that the special pencil $J=(a,b^t)R$ is a 
reduction of $I^t$. Moreover, if $R/M$ is infinite then we can take $t=1$.

\centerline{}

EXAMPLE (8.3). We can ask the question: 
if $R$ is a two dimensional normal local domain, then
what is the cardinality of $\mathfrak D(R,M)$? The answer is that it can
be any preassigned positive integer $n$. To see this let $K$ be a field and
consider the polynomial 
$$
g(X,Y,Z)=Z^m-f_1(X,Y)\dots f_n(X,Y)\in K[X,Y,Z]
$$ 
where $m>n>0$ are integers such that $m$ is nondivisible by the 
characteristic of $K$, and $f_1(X,Y),\dots,f_n(X,Y)$ are pairwise coprime 
homogeneous linear polynomials; for instance $g(X,Y,Z)=Z^3-XY$. 
Clearly $g=g(X,Y,Z)$ is irreducible in the polynomial ring $B=K[X,Y,Z]$. Let 
$$
\phi:B\to B/gB=A=K[x,y,z]
$$
be the residue class epimorphism where we have identified $K$ with $\phi(K)$ and
we have let $x=\phi(X),y=\phi(Y),z=\phi(Z)$. Upon letting $R=A_{(x,y,z)A}$ it
can easily be seen that $R$ is two dimension normal local domain with 
coefficient field $K$ and maximal ideal $M=M(R)=(x,y,z)R$.
We claim that now $\mathfrak D(R,M)$ has exactly $n$ elements.
Geometrically speaking, a section of the tangent cone of the surface $g=0$ at 
$(0,0,0)$ consists of the $n$ lines $f_1(X,Y)=0,\dots,f_n(X,Y)=0$
in the $(X,Y)$-plane passing through $(0,0)$, and they give rise to the
$n$ elements of $\mathfrak D(R,M)$.  Algebraically it can be shown that 
$\mathfrak D(R,M)$ has exactly $n$ elements thus. Let 
$$
A'=R[x',y']\;\;\text{ where }\;\;x'=x/z\;\;\text{ and }\;\;y'=y/z.
$$ 
Now
$$
z^{m-n}-f_1(x',y')\dots f_n(x',y')=0
$$
and hence, upon letting $P_i=(z,f_i(x',y'))A'$ for $1\le i\le n$, we see that
$P_1,\dots,P_n$ are exactly all the distinct height-one prime ideals in $A'$ 
which contain $MA'=zA'$ and the localizations of $A'$ at them are 
exactly all the distinct members of $\mathfrak D(R,M)$.

\centerline{}

{\it Remark $(8.3^\flat)$.} Turning to a two dimensional regular local domain,
let us cite some of the things which have been achieved.

\centerline{}

(1) In the joint paper \cite{AbH} of Abhyankar and Heinzer,
the following existence theorem of dicritical divisors is proved:
Let $R$ be a two dimensional regular local domain with
quotient field $L$. Let $U$ be any finite set of prime divisors of $R$.
Then there exists $z\in L^\times$ such that $\mathfrak D(R,z)=U$. Moreover,
if the field $R/M(R)$ is infinite then
there exists $z\in L^\times$ such that 
$\mathfrak D(R,z)^\sharp=\mathfrak D(R,z)^\flat=\mathfrak D(R,z)=U$.

\centerline{}

(2) In the joint paper \cite{AbL} of Abhyankar and Luengo, 
the following fundamental theorem of special pencils is proved:
Let $R$ be a two dimensional regular local domain with
quotient field $L$. Let $z\in L^\times$ be such that $z$ generates a special 
pencil at $R$. Then $z$ generates a polynomial pencil in $R$.

\centerline{}

(3) In Notes (3.5)(I) to (3.5)(IV) of the Abhyankar paper \cite{Ab11},
concrete examples are given to illustrate the above Remark $(8.2^\flat)$.

\centerline{}

(4) In Propositions (4.1) to (4.3) of the Abhyankar paper \cite{Ab10} and
Propositions (3.1) to (3.4) of the Abhyankar paper \cite{Ab11},
the sets
$$
\mathfrak B(R,J)^\sharp\subset\mathfrak B(R,J)^\flat
\subset\mathfrak B(R,J)\subset\mathfrak Q(R,J)\subset Q(R)
$$
of a special pencil $J$ in a two dimensional regular local domain $R$, mentioned
at the end of Section 5, are studied.

\centerline{}

{{\bf Section 9: Dicriticals of Higher Dimensional Local Domains.}} Lemmas
(6.11) to (6.14) of Section 6 together with Theorem (8.1) of Section 8
constitute the initiation of the higher dimensional dicritical theory mentioned 
in the Introduction.  Here is a scheme of how this is expected to 
be used in attacking the higher dimensional jacobian conjecture; for recent work
on this conjecture see Abhyankar's papers \cite{Ab5} to \cite{Ab7}.

Geometrically speaking, the 
possible relationship between dicritical divisors and the
jacobian conjecture which is to be exploited may be described thus.
A polynomial map from $\mathbb C^n$ to $\mathbb C^n$ is given by
$$
Y_1=f_1(X_1,\dots,X_n),\dots,Y_n=f_n(X_1,\dots,X_n)
$$
where $X_1,\dots,X_n$ are coordinates in the source $\mathbb C^n$ and
$Y_1,\dots,Y_n$ are coordinates in the target $\mathbb C^n$. For $1\le j\le n$
we have the polynomial map
$$
f_j(X_1,\dots,X_n):\mathbb C^n\to\mathbb C^1.
$$
Going over to projective spaces we get the corresponding rational map
$$
\phi_j(X_1,\dots,X_n):\mathbb P^n\to\mathbb P^1.
$$
Let $Z_1,\dots,Z_n$ be local coordinates at a point 
$\pi\in\mathbb P^n\setminus\mathbb C^n$ and let $R$ be the local ring of $\pi$
on $\mathbb P^n$. Then $R$ is an $n$-dimensional regular local domain
with maximal ideal $M=M(R)=(Z_1,\dots,Z_n)R$ and 
$$
\phi(X_1,\dots,X_n)=\frac{a_j(Z_1,\dots,Z_n)}{b_j(Z_1,\dots,Z_n)}
$$
with
$$
\text{$a_j=a_j(Z_1,\dots,Z_n)$ and $b_j=b_j(Z_1,\dots,Z_n)$ in 
$\mathbb C[Z_1,\dots,Z_n]$.}
$$
We get a pencil $I_j=(a_j,b_j)R$ in $R$ and we can consider the dicritical set
$\mathfrak D(R,I_j)$. Let us call the $n$-tuple $(f_1,\dots,f_n)$ a
jacobian $n$-tuple if the jacobian of $f_1,\dots,f_n$ with respect to
$X_1,\dots,X_n$ is a nonzero constant, and let us call it an automorphic
$n$-tuple if $\mathbb C[f_1,\dots,f_n]=\mathbb C[X_1,\dots,X_n]$.
The chain rule tells us that every automorphic $n$-tuple is a jacobian $n$-tuple.
The jacobian conjecture predicts that every jacobian $n$-tuple is an
automorphic $n$-tuple. It is plausible that 
if $(f_1,\dots,f_n)$ is a jacobian $n$-tuple then
the dicritical sets $\mathfrak D(R,I_j)_{1\le j\le n}$ are somehow related to
each other and this may help us to prove that 
$(f_1,\dots,f_n)$ is an automorphic $n$-tuple.

\end{document}